\documentclass[12pt]{article}
\usepackage{amsmath}
\usepackage{graphicx,psfrag,epsf}
\usepackage{enumerate}
\usepackage{natbib}
\usepackage{algorithm, algorithmic}
\usepackage{caption}
\usepackage[colorlinks, linkcolor=black, citecolor=blue]{hyperref}

\begin{document}
\def\spacingset#1{\renewcommand{\baselinestretch}%
{#1}\small\normalsize} \spacingset{1}
\renewcommand{\thefootnote}{\fnsymbol{footnote}}
\title{An introduction on the multivariate normal-ratio distribution}
\author{Sheng Yang\thanks{Equal contribution} \textsuperscript{1} and Zhengtao Gui$^{*}$\textsuperscript{2} \\
\textsuperscript{1}Department of Mathematics and Statistics, Lanzhou University \\
\textsuperscript{2}Department of Statistics and Finance, University of Science and Technology of China}
\date{}
\maketitle

\spacingset{1.45}

\section{Introduction}

The statistical distribution of the ratio of two normal random variables, $X$ and $Y$, is characterized by its heavy-tailed nature and absence of finite moments. The shape of its density function is highly variable, capable of exhibiting unimodal or bimodal and symmetrical or asymmetrical structures, and, in certain regions proximal to its mode, it may approximate the characteristics of a normal distribution \citep{marsaglia2006ratios,marsaglia1965ratios,diaz2013existence,hinkley1969ratio}. Several investigations have proffered diverse representations of the cumulative distribution function utilizing the framework of the bivariate normal distribution and Nicholson’s V function \citep{nicholson1943probability}. The purpose of this paper is to extend these results to the multivariate case. That is to say, suppose that the random vector $X$, where $X^{\prime}=\left(x_1, x_2, \cdots, x_p\right)$, has the multivariate normal density function:

$$
f\left(x_1, \cdots, x_p\right)=\frac{1}{(2 \pi)^{\frac{p}{2}}|\boldsymbol{\Sigma}|^{-\frac{1}{2}}} e^{-\frac{1}{2}(\boldsymbol{x}-\boldsymbol{\mu})^T \boldsymbol{\Sigma}^{-1}(\boldsymbol{x}-\boldsymbol{\mu})}
$$

We shall find the density function of the random vector $Y$, where $Y^{\prime}=\left(y_1, y_2, \cdots\right.$, $\left.y_{p-1}\right)=\left(x_2 / x_1, x_3 / x_1, \cdots, x_p / x_1\right)$. We call this distribution multivariate normal-ratio distribution.

\section{The density function in various situation}  

\label{sec:meth}

\subsection{The density function when $p$ is odd}

In our formulation, we designate $c=\frac{1}{(2 \pi)^{\frac{p}{2}}|\boldsymbol{\Sigma}|^{-\frac{1}{2}}}$, where $|\boldsymbol{\Sigma}|$ denotes the determinant of the covariance matrix $\boldsymbol{\Sigma}$. We define the transformation of variables such that $z=x_1$ and $y_i = \frac{x_{i+1}}{x_1}$ for $i=1, \ldots, p-1$. It follows that the Jacobian determinant associated with this transformation is given by $z^{p-1}$. Furthermore, we express the vector $\mathbf{X}$ as a product of $z$ and $\mathbf{W}$, specifically, $\mathbf{X} = z\mathbf{W}$, where $\mathbf{W}'=\left(1, y_1, y_2, \ldots, y_{p-1}\right)$. Consequently, the joint probability density function of the transformed variables $Y$ and $z$ is

\begin{equation}
h\left(z, Y\right)=c e^{-\frac{1}{2}(\boldsymbol{zW}-\boldsymbol{\mu})^T \boldsymbol{\Sigma}^{-1}(\boldsymbol{zW}-\boldsymbol{\mu})}z^{p-1}
\end{equation}

The density function of $Y$ is

\begin{equation}
g(Y)=c e^{\frac{1}{2}} \int_{-\infty}^{\infty} e^{-\frac{1}{2} M\left[z^{2}+(K / M) z+L / M\right]}  z^{p-1} dz
\end{equation}

where $M={W^{\prime}} \Sigma^{-1} W, K=-2 W^{\prime} \Sigma^{-1} \mu$, and $ L=u^{\prime} \Sigma^{-1} \mu+1 $.

We denote $a$ as $a=\left(L/M\right)-\left(K^{2}/4M^{2}\right)$. It can be shown that $a$ is strictly positive by examining the expression $\left(L/M\right)-4\left(K^{2}/4M^{2}\right)=\left(W^{\prime} \Sigma^{-1} W\right)^{-1}+\left(W^{\prime} \Sigma^{-1} W\right)^{-2}[\left(W^{\prime} \Sigma^{-1} W\right)\left(\mu^{\prime} \Sigma^{-1} \mu\right)
$
$
-\left(W^{\prime} \Sigma^{-1} \mu\right)^{2}] $. The term $\boldsymbol{\Sigma}^{-1}$ is positive definite, and $\left(W^{\prime} \Sigma^{-1} W\right)\left(\mu^{\prime} \Sigma^{-1} \mu\right)
-\left(W^{\prime} \Sigma^{-1} \mu\right)^{2} \geq 0$. \citep{rao1973linear}

Let $ u=(z-b) / a $, where  $b=-K /2 M$  and  $a=(L / M)-K^{2} /4 M^{2} $

\begin{equation}
g(Y)=c e^{\frac{1}{2}} a \int_{-\infty}^{\infty} e^{-\frac{1}{2} M\left(a^{2} u^{2}+a\right)}[a u+b]^{p-1} d u 
\nonumber 
\end{equation}
\begin{equation}
=c e^{\frac{1}{2}} a e^{-\frac{1}{2} M a} \int_{-\infty}^{\infty} e^{-\frac{1}{2} M a^{2} u^{2}}[a u+b]^{p-1} d u
\nonumber 
\end{equation}
\begin{equation}
=c e^{\frac{1}{2}} a e^{-\frac{1}{2} M a} \int_{-\infty}^{\infty} e^{-\frac{1}{2} M a^{2} u^{2}} \sum_{i=0}^{p-1} \binom{p-1}{i}(a u)^{p-1-i} b^{i} d u \\
\end{equation}

since $p$ is an odd integer, we have

\begin{equation}
\int_{-\infty}^{\infty} e^{-\frac{1}{2} M a^{2} u^{2}}  u^{p-2 j} d u=0 \quad \quad \text{for} \hspace{0.5em} j=1,2, \ldots,(p-1) / 2 \\
\end{equation}

In this case we only need to consider the following integral

\begin{equation}
\int_{-\infty}^{\infty} e^{-\frac{1}{2} M a^{2} u^{2}}  u^{2 j} d u\quad \quad \text{for} \hspace{0.5em} j=1,2, \ldots,(p-1) / 2 \\
\end{equation}

Then we use Euler integral to illustrate the existence of the above integral and derive its expression.

\begin{equation}
\int_{-\infty}^{\infty} e^{-t^{2}}  t^{2 j} d t=\Gamma\left(j+\frac{1}{2}\right) 
\end{equation}

\begin{equation}
\int_{-\infty}^{\infty} e^{-\frac{1}{2} M a^{2} u^{2}} u^{2 j} d u=\left(\frac{1}{2} M a^{2}\right)^{-\left(\frac{1}{2}+j\right)} \Gamma\left(j+\frac{1}{2}\right)
\end{equation}

Thus, we have that the density function of  $Y$  is 
\begin{equation}
g(Y)=c a e^{-\frac{1}{2} M a+\frac{1}{2}} \sum_{j=0}^{(p-1) / 2}\binom{p-1}{p-1-2j} a^{2 j} b^{p-1-2 j}\left(\frac{1}{2} M a^{2}\right)^{-\left(\frac{1}{2}+j\right)} \Gamma\left(j+\frac{1}{2}\right)
\end{equation}

\subsection{The density function when $p$ is even}
In the scenario where $p$ is an even integer, we employ an analogous transformation to that previously delineated. We find the Jacobian is $|z|^{p-1}$. The density of $Y$ is

\begin{equation}
g(Y)=c e^{\frac{1}{2}} \int_{-\infty}^{\infty} e^{-\frac{1}{2} M\left[z^{2}+(K / M) z+L / M\right]}  |z|^{p-1} dz
\end{equation}

Let $ u=(z-b) / a $,

\begin{equation}
g(Y)=c e^{\frac{1}{2}} a b^{p-1} \left\{ \sum_{i=0}^{p-1} \binom{p-1}{p-i-1} b^{-i} a^{i}\left[\int_{-b / a}^{\infty} e^{-\frac{1}{2} M\left(a^{2} u^{2}+a\right)}  u^{i} d u-\int_{-\infty}^{-b / a} e^{-\frac{1}{2} M\left(a^{2} u^{2}+a\right)} u^{i} d u\right]\right\}
\nonumber
\end{equation}

If $b<0$, 

\begin{equation}
\begin{split}
g(Y) = c e^{\frac{1}{2}} a b^{p-1}\left\{\sum _ { i = 0 } ^ { p - 1 }\binom{p-1}{p-i-1}b ^ { - i } a ^ { i } \left[\int_{-b / a}^{\infty} e^{-\frac{1}{2} M\left(a^{2} u^{2}+a\right)} u^{i} d u\right.\right.
\left.\left.-\int_{b / a}^{-b / a} e^{-\frac{1}{2} M\left(a^{2} u^{2}+a\right)}  u^{i} d u \right.\right. \\
\left.\left. -\int_{-\infty}^{b/a} e^{-\frac{1}{2} M\left(a^{2} u^{2}+a\right)} u^{i} d u \right]\right\} 
\end{split}
\end{equation}

For the case where $i$ is an odd integer, we have
\begin{equation}
\int_{-b/a}^{\infty} e^{-\frac{1}{2} M\left(a^{2} u^{2}+a\right)} u^{i} d u -\int_{-\infty}^{b/a} e^{-\frac{1}{2} M\left(a^{2} u^{2}+a\right)} u^{i} d u = 2 \int_{-b / a}^{\infty} e^{-\frac{1}{2} M\left(a^{2} u^{2}+a\right)} u^{i} d u
\nonumber
\end{equation}
\begin{equation}
\int_{b / a}^{-b / a} e^{-\frac{1}{2} M\left(a^{2} u^{2}+a\right)} u^{i} d u = 0
\nonumber
\end{equation}

For the case where $i$ is an even integer, we have

\begin{equation}
\int_{-b/a}^{\infty} e^{-\frac{1}{2} M\left(a^{2} u^{2}+a\right)} u^{i} d u -\int_{-\infty}^{b/a} e^{-\frac{1}{2} M\left(a^{2} u^{2}+a\right)} u^{i} d u=0 
\nonumber
\end{equation}
\begin{equation}
\int_{b / a}^{-b / a} e^{-\frac{1}{2} M\left(a^{2} u^{2}+a\right)} u^{i} d u = 2\int_{0}^{-b / a} e^{-\frac{1}{2} M\left(a^{2} u^{2}+a\right)} u^{i} d u
\nonumber
\end{equation}

Accordingly, it follows that the density function of 
$Y$ can be delineated as

\begin{equation}
\begin{split}
g(Y)=2 c e^{\frac{1}{2}} a b^{p-1} e^{-\frac{1}{2} M a}\left\{\sum_{i=1}^{p / 2}\binom{p-1}{p-2i} a^{2 i-1} b^{-(2 i-1)} \int_{-b / a}^{\infty} e^{-\frac{1}{2} M a^{2} u^{2}} u^{2 i-1} d u\right. \\
\left.+\sum_{i=0}^{(p-2) / 2}\binom{p-1}{p-2i-1}
a^{2 i} b^{-2 i} \int_{0}^{-b / a} e^{-\frac{1}{2} M a^{2} u^{2}} u^{2 i} d u\right\}
\end{split}
\end{equation}

In the event that $b>0$, one may employ an analogous approach to deduce the density function of $Y$

\begin{equation}
\begin{split}
g(Y)=2 c e^{\frac{1}{2}} a b^{p-1} e^{-\frac{1}{2} M a}\left\{\sum_{i=1}^{p / 2}\binom{p-1}{p-2i} a^{2 i-1} b^{-(2 i-1)} \int_{b / a}^{\infty} e^{-\frac{1}{2} M a^{2} u^{2}} u^{2 i-1} d u\right. \\
\left.+\sum_{i=0}^{(p-2) / 2}\binom{p-1}{p-2i-1}
a^{2 i} b^{-2 i} \int_{0}^{b / a} e^{-\frac{1}{2} M a^{2} u^{2}} u^{2 i} d u\right\}
\end{split}
\end{equation}

Subsequently, we address the evaluation of the infinite and definite integrals presented in the preceding
expression.

Under the condition that $b>0$, it follows that 

\begin{equation}
 \int_{b / a}^{\infty} e^{-\frac{1}{2} M a^{2} u^{2}} u^{2 i-1} d u
 = \int_{0}^{\infty} e^{-\frac{1}{2} M a^{2} u^{2}} u^{2 i-1} d u-
 \int_{0}^{b / a} e^{-\frac{1}{2} M a^{2} u^{2}} u^{2 i-1} d u
\end{equation}

\begin{equation}
\int_{0}^{\infty} e^{-\frac{1}{2} M a^{2} u^{2}} u^{2 i-1} d u=\left(\frac{1}{2} M a^{2}\right)^{-i} \frac{1}{2} \Gamma(i)
 \end{equation}

For the case where $i=0$,  the integral 
$\int_{0}^{b/a} e^{-\frac{1}{2} M a^{2} u^{2}} u^{2 i} d u$ simplifies to $\left(\frac{1}{2} M a^{2}\right)^{-\frac{1}{2}} \int_{0}^{d} \cdot e^{-t^{2}} d t $, where $d$ is defined as $\left(\frac{1}{2} M a^{2}\right)^{\frac{1}{2}} \frac{b}{a}$ .

For $i \neq 0 $, the integral $\quad \int_{0}^{b/a} e^{-\frac{1}{2} M a^{2} u^{2}} u^{2 i} d u$ can be expressed as $\left(\frac{1}{2} M a^{2}\right)^{-\left(\frac{1}{2}+i\right)} \frac{1}{2} \int_{0}^{c} e^{-z} z^{i-\frac{1}{2}} d z $, $i=1, \ldots(p-2) / 2$
, where $c$ is defined as $d^{2}$ .
In the case at hand, the complexity of the integrals can be ameliorated through the strategic application of integration by parts.

\section{A simple approximation}

Consider the random vector $X$, where $X^{\prime}=\left(x_1, x_2, \cdots, x_p\right)$ follows a multivariate normal distribution with $\mu=\left(\mu_1, \mu_2, \cdots, \mu_p\right)$ and $\Sigma$. We present a succinct condition under which the distribution of a vector $Y$, where $Y^{\prime}=\left(y_1, y_2, \cdots\right.$, $\left. y_{p-1}\right)=\left(x_2 / x_1, x_3 / x_1, \cdots, x_p / x_1\right)$ closely approximates a multivariate normal distribution.

Assuming that the variable \(x_1\) remains invariably and strictly greater than zero, it follows that the set of problems under consideration can be resolved perfectly since $\Pr(Y < \boldsymbol{t}) = \Pr\left(y_1 < t_1, y_2 < t_2, \cdots, y_{p-1} < t_{p-1}\right) = \Pr\left(x_2 - x_1 t_1 < 0, x_3 - x_1 t_2 < 0, \cdots, x_p - x_1 t_{p-1} < 0\right)$.

Given \( x_1, x_2, \cdots, x_p \) are normal random variables, and \( x_2 - x_1 t_1, x_3 - x_1 t_2, \cdots, x_p - x_1 t_{p-1} \) are linear combinations of them, they can be expressed as a random variable of a normal distribution. If we know the mean, variance, and covariance of the random variables \( x_2 - x_1 t_1, x_3 - x_1 t_2, \cdots, x_p - x_1 t_{p-1} \), we can use the cumulative distribution function of the multivariate normal distribution to obtain $\Pr(Y < \boldsymbol{t})$.

For example, given three normal random variables \(x\), \(y\), \(z\) having $\mu_x = 10$, $\mu_y = 0$, $\mu_x = 0$, unit variance (\(\sigma^2 = 1\)), and a pairwise correlation coefficient \(\rho_{xy} = \rho_{xz} = \rho_{yz} = 0.5\). We can consider that the distribution of the random vector $Y$, where $Y^{\prime} = \left(y/x, z/x\right)$ is very close to the bivariate normal distribution since the probability of $x$ being less than 0 is small.

Assuming $\boldsymbol{t} = (2,3)$, we have $\Pr(Y < \boldsymbol{t}) \approx \Pr\left(y - 2x < 0, z - 3x < 0\right)$ and 

\[ E[y - 2x] = E[y] - 2E[x] = -20 \]
\[ E[z - 3x] = E[z] - 3E[x] = -30 \]
\[ \text{Var}(y - 2x) = \text{Var}(y) + 4\text{Var}(x) - 4\text{Cov}(x, y) = 3 \]
\[ \text{Var}(z - 3x) = \text{Var}(z) + 9\text{Var}(x) - 6\text{Cov}(x, z) = 7 \]
\[ \text{Cov}(y - 2x, z - 3x) = \text{Cov}(y, z) - 2\text{Cov}(x, z) - 3\text{Cov}(x, y) + 6\text{Var}(x) = 4.5 \]

This is a bivariate normal distribution with $\mu = (-20, -30)$ and $\Sigma =
\begin{bmatrix}
3 & 4.5 \\
4.5 & 7
\end{bmatrix}$. We can use the cumulative distribution function (CDF) of the bivariate normal distribution \citep{owen1956tables} to calculate the joint probability $\Pr\left(y - 2x < 0, z - 3x < 0\right)$.

\section{Applications}

Consider an economic model where in the demand for a commodity, denoted by \( y \), is a linear function of its price \( p \) and income \( m \). Specifically, this relationship can be represented as the following linear regression model:

\begin{equation}
y = \beta_0 + \beta_1 p + \beta_2 m + \epsilon
\nonumber
\end{equation}

where \( \epsilon \) denotes the error term, assumed to be identically and independently distributed following a normal distribution: \( \epsilon \sim N(0, \sigma^2) \).

The key interest lies in assessing the sensitivity of the commodity demand to changes in price and income, known in economics as elasticities. Price elasticity and income elasticity are defined as:

\begin{equation}
\varepsilon_p = \frac{\partial y / y}{\partial p / p} = \frac{\beta_1 p}{y}
\nonumber
\end{equation}

\begin{equation}
\varepsilon_m = \frac{\partial y / y}{\partial m / m} = \frac{\beta_2 m}{y}
\nonumber
\end{equation}

Given that both \( y \) and \( p \) are considered random variables and are normally distributed, the distribution of the elasticity \( \varepsilon_p \) is uncertain and is no longer normal.

To delve further into this, assume that the price \( p \), demand \( y \) and income \( m \) follow normal distributions: \( p \sim N(\mu_p, \sigma_p^2) \) and \( y \sim N(\mu_y, \sigma_y^2) \) and \( m \sim N(\mu_m, \sigma_m^2) \). Under this assumption, the joint distribution of the \( \varepsilon_p \) and \( \varepsilon_m \) is the multivariate normal-ratio distribution. Determining the distribution of the random vector (\( \varepsilon_p \), \( \varepsilon_m \)) requires a deep investigation into the properties of this distribution. 

In conclusion, the joint distributional properties of price elasticity and income elasticity might be intricate. Grasping the characteristics of multivariate normal-ratio distribution is pivotal to guaranteeing the precision of economic analyses and making informed decisions. This distribution offers a methodology to investigate and comprehend the complexities of these ratios in high-dimensional spaces.

\section{Discussion}

In this article, we study the multivariate normal ratio distribution and propose a concise condition under which this distribution will approximate the multivariate normal distribution. The exploration of the asymptotic and moment properties of the multivariate normal ratio distribution has always been a challenging topic in the future. Marsaglia's study \citep{marsaglia1965ratios} evaluates the density and distribution functions of the ratio $z/w$ for any two jointly normal variates $z, w,$ and provides details on methods for transforming a general $z/w$ into a standard form $(a+x)/(b+y)$, with $x$ and $y$ independent standard normal and $a, b$ non-negative constants. It discusses handling general ratios when, in theory, none of the moments exist yet practical considerations suggest there should be approximations. These approximations show that many of the ratios of normal variates encountered in practice can themselves be taken as normally distributed. Whether the multivariate normal ratio distribution can converge to the multivariate normal distribution under the constraints of other conditions still needs to be explored by researchers.

\bibliographystyle{elsarticle-harv}
 
\bibliography{mybib}

\end{document}